\documentclass[11pt,leqno]{article}

\usepackage{amssymb,amsmath,amsthm}

\theoremstyle{plain}
\newtheorem{thm}{Theorem}

\newtheorem{lem}[thm]{Lemma}

\theoremstyle{remark}
\newtheorem*{rk}{Remark}

\theoremstyle{definition}

\newcommand{\ovl}{\overline}
\newcommand{\bb}[1]{\mathbb{#1}}
\newcommand{\cl}[1]{\mathcal{#1}}
\newcommand{\vp}{\varepsilon}
\newcommand{\sst}{\scriptstyle}

\begin{document}

\title{Remarks on $B(H)\otimes B(H)$}

\author{by\\ \phantom{blank line}\\ Gilles Pisier\thanks{Supported in part by 
NSF grant 
No.~0200690}\\ Department of Mathematics\\ 
Texas A\&M University\\ College Station, TX 77843-3368, 
U. S. A.\\ and\\ Universit\'e Paris VI\\ Equipe d'Analyse, Case 186, 75252\\ 
Paris Cedex 05, France}

\date{}
\maketitle

\begin{abstract}
We review the existing proofs that the min and max norms are different 
on $B(H)\otimes B(H)$ and give a shortcut avoiding the consideration of 
non-separable families of operator spaces.
\end{abstract}

Math. Subject. Class.: 46L05, 46M05, 47D15
\bigskip

Let $A,B$ be $C^*$-algebras. Let $\|~~\|_\alpha$ be a $C^*$-norm on their 
algebraic tensor product, denoted by $A\otimes B$; as usual, $A\otimes_\alpha B$ 
then denotes the $C^*$-algebra obtained by completing $A\otimes B$ with respect 
to $\|~~\|_\alpha$. By classical results (see \cite{T}) the set of $C^*$-norms 
admits a minimal and a maximal element denoted respectively by 
$\|\cdot\|_{\min}$ and $\|\cdot\|_{\max}$. Then $A$ is called nuclear if for any 
$B$ we have $A \otimes_{\min} B = A \otimes_{\max} B$, or equivalently 
$\|x\|_{\min} = \|x\|_{\max}$ for any $x$ in $A\otimes B$. We refer the reader 
to \cite{T,KR} for more information on nuclear $C^*$-algebras. We note in 
particular that by results due to Connes and Haagerup (\cite{C,H}), a 
$C^*$-algebra is nuclear iff it is amenable as a Banach algebra (in 
B.E.~Johnson's sense).

Although, it was known early (\cite{W}) that $B(H)$ is not nuclear (assuming 
$\dim(H)=\infty$), the problem whether $A \otimes_{\min} B = A \otimes_{\max} B$ 
when $A=B=B(H)$ remained open until \cite{JP}. In the latter paper, several 
different proofs were given. Taking into account the most recent information 
from \cite{HT}, we now know that
\begin{equation}\label{eq1}
\sup\left\{\frac{\|t\|_{\max}}{\|t\|_{\min}} ~\Big|~ t\in B(H) \otimes B(H) 
\quad {\rm rk}(t)\le n\right\} \ge
\frac{n}{2\sqrt{n-1}}.
\end{equation}
This estimate is rather sharp asymptotically, since it can be shown that the 
supremum appearing in \eqref{eq1} is $\le \sqrt n$ (see \cite{JP} or 
\cite[p.~353]{P1}). 

\begin{rk}
In \eqref{eq2} below, the complex conjugate $\bar a$ of a matrix $a$ in $M_N$ is 
meant in the usual way, i.e.\ $(\bar a)_{ij} = \ovl{a_{ij}}$.  In general, we 
will need to consider the conjugate $\bar A$ of a $C^*$-algebra $A$. This is the 
same object but with the complex multiplication changed to $(\lambda,a)\to 
\bar\lambda a$, so that $\bar A$ is anti-isomorphic to $A$.
For any $a\in A$, we denote by $\bar a$ the same element
 viewed as an element of $\bar A$.
 Note that $\bar A$ 
can also be identified with the opposite $C^*$-algebra $A^{op}$ which is defined 
as the same object but with the product changed to $(a,b)\to ba$. It is easy to 
check that the mapping $a\to a^*$ 
(or more properly $\bar a\to a^*$) is a (linear) $*$-isomorphism from $\bar A$ to 
$A^{op}$.

The distinction between $A$ and $\bar A$ is necessary in general, but not for $A 
= B(H)$ since in that case, using $H\simeq \ovl H$,
we have $\ovl{B(H)} \simeq  B(\ovl H) \simeq
B(H)$, and in particular  $\ovl{M_N}  \simeq
M_N $. 
As a consequence, 
it is easy to see that for any matrix $a$ 
in $M_N(A)$ we have $$\|[\bar a_{ij}]\|_{M_N(\bar A)}
 = \|[a_{ij}]\|_{M_N(A)}.$$

Note however that $H\simeq \ovl H$ depends
on the choice  of a basis so the isomorphism
$\ovl{B(H)}\simeq B(H)$ is not canonical.
Nevertheless, this shows that the problem whether the min and max norms are the 
same is identical for $B(H) \otimes B(H)$ and for $B(H) \otimes \ovl{B(H)}$.
\end{rk}

\begin{rk}
Consider $a_1,\ldots, a_n$ in $A$ and $b_1,\ldots, b_n$ in $B$. Using the 
preceding remark we have
\[
\left\|\sum a_i\otimes\bar b_i\right\|_{A\otimes_\alpha\ovl B} = \left\|\sum a_i 
\otimes b^*_i\right\|_{A\otimes_\alpha B^{op}}
\]
for any $C^*$-norm, in particular for $\alpha = \min$ or max. Moreover, we have
\[
\left\|\sum a_i\otimes \bar b_i \right\|_{A\otimes_{\max}\ovl B} = \left\|\sum 
a_i \otimes b^*_i\right\|_{A\otimes_{\max} B^{op}} = \sup\left\{\left\|\sum 
\pi(a_i) \sigma(b_i)^*\right\|\right\}
\]
where the supremum runs over all commuting range pairs $\pi\colon \ A\to B(H)$, 
$\sigma\colon \ B\to B(H)$ with $\pi$ a representation and $\sigma$ an 
anti-representation on the same (arbitrary) Hilbert space $H$. 
\end{rk}

\begin{rk}
Let $M$ be a $C^*$-algebra equipped with a tracial state $\tau$. Then the GNS 
construction associated to $(M,\tau)$ produces a Hilbert space $H$, a cyclic 
unit vector $\xi$ in $H$ and commuting left and right actions of $M$ on $H$ (we 
denote the latter simply by $a\cdot h\cdot b$ for $h\in H$, $a,b\in A$) so that 
$x\cdot \xi 
= \xi\cdot x$ and $\tau(xy) = \langle x\cdot \xi\cdot y,\xi\rangle$ $(x,y\in 
A)$. If we denote $L(a)h = a\cdot h$ and $R(a)h=h\cdot a$ then $L$ (resp.\ $R$) 
is a representation of $A$ resp.\ $(A^{op})$ on $B(H)$ and the ranges of $L$ and 
$R$ commute.

We have then for any $n$-tuple $(u_1,\ldots, u_n)$ in $M$
\begin{equation}\label{eq0}
\left\|\sum\nolimits^n_1 u_i \otimes \bar u_i\right\|_{M\otimes_{\max} \ovl M} = 
\left\|\sum\nolimits^n_1 u_i\otimes u^*_i\right\|_{M \otimes_{\max} M^{op}} = 
n. \end{equation}
Indeed, this is
\begin{align*}
\ge \left\|\sum\nolimits^n_1 L(u_i)R(u^*_i)\right\| &\ge 
\left\|\sum\nolimits^n_1 u_i\cdot\xi\cdot u^*_i\right\|\\
\end{align*}
but $  u_i\cdot\xi\cdot u^*_i = \xi $ hence 
\eqref{eq0} is $\ge n$ and $\le n$ is trivial by the triangle inequality.
 In particular, for any unitary matrices $u_1,\ldots, u_n$ in $M_N$ we have
\[ \left\|\sum\nolimits^n_1 u_i\otimes \bar u_i\right\|_{\max} = n.\]

\end{rk}

In all the proofs in \cite{JP}, a crucial role is played by a certain constant 
$C(n)$, defined as follows:\ $C(n)$ is the smallest constant $C$ such that for 
each $m\ge 1$, there is $N_m\ge 1$ and an $n$-tuple $[u_1(m), \ldots, u_n(m)]$ 
of unitary $N_m\times N_m$ matrices such that
\begin{equation}\label{eq2}
\sup_{m\ne m'} \left\|\sum^n_{i=1} u_i(m) \otimes\ovl{u_i(m')}\right\|_{\min} 
\le C.
\end{equation}
The crucial fact to show that $B(H) \otimes_{\min} B(H)\ne B(H) \otimes_{\max} 
B(H)$ is that $C(n) < n$. Various improvements were given over the initial 
estimates of $C(n)$ in \cite{JP}, notably by A. Valette 
using Ramanujan graphs (see \cite{P1} for more
details) . The final word on this is now

\begin{thm}[\cite{HT}]\label{thm1}
$C(n) = 2\sqrt{n-1}$ for any $n\ge 2$. \end{thm}

The (much easier) lower bound 
$2\sqrt{n-1} \le C(n)$ was proved in \cite{P2}.

The connection of $C(n)$ to $B(H) \otimes B(H)$ goes through the next statement.

\begin{thm}[\cite{JP}]\label{thm2}
For any $n\ge 1$ and any $C>C(n)$, there is a tensor $t$ of rank
$n$ in $B(H) \otimes B(H)$ such  that
\[
\|t\|_{\max}/\|t\|_{\min} \ge n/C.
\]\end{thm}

Our goal is to describe below a shortcut on the proof of this last result. Our 
``shortcut'' avoids any reference to the non-separability of the space of 
$n$-dimensional operator space (as was done in \cite{JP}) and uses instead a 
compactness argument for ``convergence in distribution'' of $n$-tuples of 
operators.

More precisely, the notion of ``distribution'' that we will use is the same as 
Voiculescu's definition in free probability but our terminology is slightly 
different. Let $S$ be the set consisting of the disjoint union of the sets 
\[
S_k = [1,\ldots, n]^k \times\{1,*\}^k.
\]
For any $w = \left((i_1,\ldots, i_k), (\vp_1,\ldots, \vp_k)\right)$ in $S_k$ 
and any $n$-tuple $x =(x_1,\ldots, x_n)$ in $B(H)$ we denote
\[
w(x) = x^{\vp_1}_{i_1} x^{\vp_2}_{i_2} \ldots
x^{\vp_k}_{i_k}
\]
(where $x^\vp = x$ if $\vp=1$ and $x^\vp= x^*$ if $\vp=*$). Let $x = 
(x_1,\ldots, 
x_n)$ be an $n$-tuple in a von~Neumann algebra $M$ equipped with a tracial state 
$\tau$. By ``the distribution of $x$,'' we mean the function
\[
\mu_x\colon \ S\to {\bb C}
\]
defined by
\[
\mu_x(w) = \tau(w(x)).
\]
When $x=(u_1,\ldots, x_n)$ is an $n$-tuple of unitary operators, we may as well 
consider that $\mu_x$ is a function defined on ${\bb F}_n$ (free group with 
generators $g_1,\ldots, g_n$) by setting $\mu_x(w) = \tau(\pi_x(w))$ where 
$\pi_x\colon \ {\bb F}_n\to M$ is the unitary representation defined by 
$\pi(g_i) = u_i$.

The following is elementary and well known.

\begin{lem}\label{lem3}
Fix $n$. Let $(M(m), \tau_m)$ be a sequence of von~Neumann algebras equipped 
with (tracial) states. Let $x(m) = (x_1(m),\ldots, x_n(m))$ be a bounded 
sequence of $n$-tuples with $x(m)\in M(m)^n$. Then there is a subsequence 
$\{m_k\}$ such that the distributions of $x(m_k)$ converge pointwise on $S$ when 
$k\to \infty$.
\end{lem}

To identify the limit of a sequence of distribution, it will be convenient to 
use von~Neumann algebra ultraproducts. We briefly recall how those are 
constructed.

Let $(M(m),\tau_m)$ be as before $(m\in{\bb N})$. Let ${\cl U}$ be a non-trivial 
ultrafilter on ${\bb N}$. Let ${\cal B}= \bigoplus\limits_m M(m)$. We equip the 
von~Neumann algebra ${\cal B}$ with the state $\tau_{\cl U}$ defined by
\[
\forall y = (y_m)\in {\cal B}\qquad \tau_{\cl U}(y) = \lim\nolimits_{\cl U} 
\tau_m(y_m).
\]
The GNS-construction applied to ${\cal B}$ and the state $\tau_{\cl U}$ produces a 
Hilbert space $H_{\cl U}$ and left and right actions of ${\cal B}$ on $H_{\cl U}$, each 
with kernel equal to
\[
I_{\cl U} = \{y = (y_m)\mid \lim\nolimits_{\cl U} \tau_m(y^*_my_m) = 0\}.
\]
We then set
\[
M({\cl U}) = {\cal B}/I_{\cl U}.
\]
Thus, after passing to the quotient by $I_{\cl U}$ we obtain an isometric 
representation:
\[
a \to L(a) \in B(H_{\cl U})
\]
of $M({\cl U})$ on $H_{\cl U}$ and an isometric representation $a\to R(a) \in 
B(H_{\cl U})$ of $M({\cl U})^{op}$.

It is well known (cf.\ e.g.\ \cite[p.~211]{P1}) that $M({\cl U})$ is a 
von~Neumann algebra, that $L(M({\cl U}))$ is a von~Neumann subalgebra of 
$B(H_{\cl U})$, and that we have
\[
L(M({\cl U}))' = R(M({\cl U}))\qquad R(M({\cl U}))' = L(M({\cl U})).
\]
Let $x = \{x(m)\mid m\in {\bb N}\}$ be a bounded sequence of $n$-tuples with 
$x(m)\in M(m)^n$ as before. Equivalently, $x$ can be viewed as
an $n$-tuple of elements of ${\cl B}$ (i.e. as an
element of ${\cl B}^n$). Let
$\hat  x = (\hat x_1,\ldots, \hat x_n)$ be the associated $n$-tuple in
$({\cl B}/I_{\cl  U})^n$. Then, for any ``word'' $w$ in $S$, we clearly have
\[
\tau_{\cl U}(w(\hat x)) = \lim\nolimits_{\cl U} \tau_m(w(x(m))).
\]
Hence the distribution of $x(m)$ tends pointwise to that of $\hat x$ along ${\cl 
U}$, so we can write $\lim_{\cl U} \mu_{x(m)} = \mu_{\hat x}$. The next (again 
elementary and well known) lemma connects limits in distribution with 
ultrafilters.

\begin{lem}\label{lem4}
Let $\{(x(m)\mid m\in {\bb N}\}$ be a sequence of $n$-tuples as in the preceding 
Lemma. The following are equivalent.
\begin{itemize}
\item[(i)] The distributions of $x(m)$ converge pointwise when $m\to \infty$.
\item[(ii)] For any non-trivial ultrafilter ${\cl U}$ on ${\bb N}$, the 
associated $n$-tuple $\hat x = (\hat x_1,\ldots, \hat x_n)$ in $(M({\cl U}), 
\tau_{\cl U})$ has the same distribution (i.e.\ its distribution does not depend 
on ${\cl U}$).
\item[(iii)] There is a von Neumann algebra $(M,\tau)$ equipped with a tracial 
state and $y = (y_1,\ldots, y_n)$ in $M^n$ such that $x(m)\to y$ in 
distribution.
\end{itemize}
\end{lem}

\begin{thm}\label{thm5}
Let $\{[u_1(m),\ldots, u_n(m)], m\in {\bb N}\}$ a sequence of $n$-tuples of 
unitary matrices satisfying \eqref{eq2} (recall $u_1(m),\ldots, u_n(m)$ are of 
size $N_m\times {\bb N}_m$). Suppose that $[u_1(m),\ldots, u_n(m)]$ converges in 
distribution when $m\to\infty$. Let ${\bb N} = \alpha\cup\beta$ be any disjoint 
partition of ${\bb N}$ into two infinite subsets, let $u_i(\alpha) = 
\bigoplus\limits_{m\in\alpha} u_i(m)$, $u_i(\beta) = \bigoplus\limits_{m'\in 
\beta} u_i(m')$, and finally let
\[
t =\sum^n_{i=1} u_i(\alpha)\otimes \ovl{u_i(\beta)}.
\]
We have then
\[
\|t\|_{\min} \le C\quad\text{and}\quad \|t\|_{\max} = n,
\]
hence $\|t\|_{\max}/\|t\|_{\min} \ge n/C$, where the min and max norms are 
relative to $\left(\bigoplus\limits_{m\in\alpha} M_{N_m}\right) \otimes 
\overline{ \left(\bigoplus\limits_{m\in \beta} M_{N_m}\right)}$.
\end{thm}

\begin{proof}
We have obviously
\[
\|t\|_{\min} = \sup_{\underset{\sst m'\in\beta}{m\in\alpha}} \left\|\sum u_i(m) 
\otimes \ovl{u_i(m')}\right\|
\]
hence $\|t\|_{\min}\le C$. We now turn to $\|t\|_{\max}$. Let ${\cl U}$ be a 
nontrivial ultrafilter on $\alpha$ and let ${\cl V}$ be one on $\beta$. We 
construct the ultraproducts $M({\cl U})$ and $M({\cl V})$ as above. Since the 
quotient mappings $\bigoplus\limits_{m\in\alpha} M(N_m)\to M({\cl U})$ and 
$\bigoplus\limits_{m\in\beta} M(N_m)\to M({\cl V})$ are $*$-homomorphisms, we 
have
\[
\|t\|_{\max} \ge \left\|\sum u_i\otimes \bar v_i\right\|_{M({\cl U}) 
\otimes_{\max} \overline{M({\cl V})}}
\]
where $u_i$ (resp.\ $v_i$) is the equivalence class modulo ${\cl U}$ (resp.\ 
${\cl V}$) of $\bigoplus\limits_{m\in\alpha} u_i(m)$ (resp.\ 
$\bigoplus\limits_{m\in\beta} u_i(m)$).

Now, since we assume that $[u_1(m),\ldots, u_n(m)]$ converges in distribution, 
$(u_1,\ldots, u_n)$ and $(v_1,\ldots, v_n)$ must have the same distribution 
relative respectively to $\tau_{\cl U}$ and $\tau_{\cl V}$. But this implies 
that there is a $*$-isomorphism $\pi$ from the von~Neumann algebra $M_u$ 
generated by $(u_1,\ldots, u_n)$ to the one $M_v$ generated by $(v_1,\ldots, 
v_n)$, defined simply by $\pi(u_i) = v_i$. Moreover, since we are dealing here 
with {\em finite\/} traces, there is a conditional expectation $P$ from $M({\cl 
U})$ onto $M_u$. Therefore the composition $T = \pi P$ is a  unital completely 
positive map such that $T(u_i) = v_i$. Hence we have 
\begin{align*}
\left\|\sum u_i \otimes \bar v_i\right\|_{\max} &\ge \left\|\sum  T(u_i) 
\otimes \bar v_i\right\|_{M(v) \otimes_{\max}\ovl{M(v)}}\\
&= \left\|\sum v_i \otimes \bar v_i\right\|_{M(v) \otimes_{\max} \ovl{M(v)}}.
\end{align*}
But then by the Remark preceding Theorem \ref{thm1} we conclude that $\|t\|_{\max} \ge n$.
\end{proof}

 \noindent {\it Proof of Therem \ref{thm2}}. Fix any number $C> C(n)$. Then
there is a sequence $\{[u_1(m),\ldots, u_n(m)], m\in {\bb N}\}$
 satisfying
\eqref{eq2}. By Lemma \ref{lem3}, there is a subsequence that converges in
distribution. Let $\{N(m)\}$ be the sequence of sizes for this subsequence
and let $\cl B= \oplus_m M_{N(m)}$.
Applying Theorem \ref{thm5} to this subsequence
we find $t$ in $ {\cl B}\otimes \bar{\cl B} $ such that 
\[
\|t\|_{\max}/\|t\|_{\min} \ge n/C.
\]
But since $ {\cl B}\subset B(H)$ and there is  a unital completely positive
projection from $B(H)$ onto $ {\cl B}$, the min and max norms of $t$
viewed as sitting in $  {  B(H)}\otimes \bar{  B(H)} $ are the same
as when computed in $ {\cl B}\otimes \bar{\cl B} $
(this essentially goes back to \cite{L}. see also \cite[chapter 14]{P1}).
Thus we obtain  Theorem
\ref{thm2}.

\begin{rk} The same shortcut applies to the proof of the main result in
\cite{O}.

\end{rk}

\end{document}